\documentclass[12pt]{article}
\usepackage{amsmath}
\usepackage{amsthm}
\usepackage{amsfonts}
\usepackage{amssymb}

\let\phi\varphi
\def\C{\mathbb C}

\let\phi\varphi
\let\epsilon\varepsilon
\def\C{\mathbb C}
\def\Z{\mathbb Z}
\def\N{\mathbb N}
\def\mmod{\ \mathrm{mod}\ }
\let\wt\widetilde
\def\Im{\mathop{\mathrm{Im}}}

\makeatletter

\newenvironment{proofd}[1][\proofname]{\par
  \pushQED{\qed}%
  \normalfont \topsep\smallskipamount\relax
  \trivlist
  \item[\indent\hskip\labelsep
        \itshape
    #1]\ignorespaces
}{%
  \popQED\endtrivlist\@endpefalse
}

\makeatother
\newtheorem{statement}{Proposition}

\theoremstyle{remark}

\newtheorem*{remark}{Remark}
\newtheorem*{remarks}{Remarks}

\theoremstyle{definition}

\newtheorem*{definition}{Definition}

\title{Belavin Elliptic $R$-Matrices and Exchange Algebras}
\author{A. V. Odesskii\thanks{
      This work was supported by CRDF grant RP1-2254, RFBR grants 99-01-01169 and 00-15-96579, and INTAS
grant 00-00055.
}}
\date{}

\begin{document}

\sloppy

\maketitle

\begin{abstract}

 We study Zamolodchikov algebras whose commutation relations are described by
    Belavin matrices defining a solution of the Yang--Baxter equation (Belavin $R$-matrices). Homomorphisms of Zamolodchikov algebras into dynamical algebras with exchange relations and also of
    algebras with exchange relations into Zamolodchikov algebras are constructed. It turns out that
    the structure of these algebras with exchange relations depends substantially on the primitive $n$th
    root of unity entering the definition of Belavin $R$-matrices.
\end{abstract}

\section*{Introduction}
    In this paper, we study Zamolodchikov algebras for Belavin elliptic $R$-matrices. This is a family
$Z_{n,k}(\Gamma,\eta)$ of associative algebras, where $\Gamma\subset\C$    is the integer lattice generated by the elements 1
and~$\tau$, $\Im\tau>0$, $\eta\in\C$, and $n$ and $k$ are coprime
positive integers such that $1\le k<n$. The
algebra $Z_{n,k}(\Gamma,\eta)$ is determined by the generators
$\bigl\{x_i(u);\,i\in\Z/n\Z,\,u\in\C\bigr\}$ and the relations
\begin{multline*}
\frac{\theta_1(0)\dots\theta_{n-1}(0)\theta_0(v-u+\eta)\dots
\theta_{n-1}(v-u+\eta)}{\theta_0(\eta)\dots\theta_{n-1}(\eta)
\theta_0(v-u)\dots\theta_{n-1}(v-u)}x_\alpha(u)x_\beta(v)\\
=
\sum_{ r\in\Z/n\Z}\frac{\theta_{\beta-\alpha+ r(k-1)}(v-u+\eta)}
{\theta_{k r}(\eta)\theta_{\beta-\alpha- r}(v-u)}
x_{\beta-r}(v)x_{\alpha+r}(u).
\end{multline*}
Here $\theta_\alpha(u)\in\Theta_{n,\frac{n-1}2}(\Gamma)$ is a theta function of order $n$ (see Sec.~1.1). These formulas are distinct
from the Belavin formulas [2], but it can readily be verified that the corresponding Zamolodchikov
algebras are isomorphic. In the Belavin notation, the role of $k$ is played by a primitive $n$th root
of unity.

    As is known, in studying integrable models related to elliptic $R$-matrices, there appear some additional difficulties (in comparison with the trigonometric and rational cases). The same difficulties
are encountered in studying the representations of the algebra
$Z_{n,k}(\Gamma,\eta)$ or of the corresponding
algebra of $L$-operators.

    The fact is that these algebras have no analog of Cartan subalgebras, and therefore it is unclear
how to construct an analog of highest weight modules since there is no notion of weight at all.

    A method for overcoming these difficulties was found by Baxter [1] in studying $XYZ$-models
(the case $n = 2$, $k = 1$). In the language of [1], this is the reduction of the $XYZ$-model to the
$SOS$-model (the so-called vertex-face correspondence). In the
algebraic language, this is the construction of homomorphisms from the
algebra $Z_{n,k}(\Gamma,\eta)$ into some other algebras for which highest
weight modules already exist. For $k = 1$, these are Zamolodchikov algebras for dynamical elliptic
$R$-matrices (see [4]). For the other values of $k$, these algebras have a more complicated structure, and
no related models are probably known. Let us first describe these algebras for the case $k = 1$. These
are the algebras $X_1^m(\Gamma,\mu;\lambda)$ generated by the commutative subalgebra consisting of meromorphic
functions of the variables $y_1, \dots , y_m$ and by the generators $e_1(u), \dots , e_m(u)$, $u\in   \C$. The relations
have the form $e_\alpha(u)f(y_1,\dots,y_m)=f(y_1+\lambda,\dots,
y_\alpha+\lambda-\mu,\dots,y_m+\lambda)e_\alpha(u)$, i.e., $y_1, \dots , y_m$ are
dynamical variables. Here  $\lambda,\mu\in\C$ are fixed. Moreover,
$e_\alpha(u)e_\beta(v)=
\phi_1e_\beta(v)e_\alpha(u)+\phi_2e_\alpha(v)e_\beta(u)$,
where $\alpha\ne\beta$, and $\phi_ 1$ and  $\phi_2$ are meromorphic functions of the variables $u, v, y_\alpha$, and $y_\beta$  (ratios
of some theta functions). If    $\alpha=\beta$, then $e_\alpha(u)e_\alpha(v)=e_\alpha(v)e_\alpha(u)$. For each $m\in\N$, there is a
homomorphism $x_i(u)\mapsto\sum_{1\le\alpha\le m}\theta_i(y_\alpha+u)e_\alpha(nu)$ from $Z_{n,1}(\Gamma,\eta)$ into an algebra of the above
type. Here $\theta_i(z)$ stands for a theta function of order $n$ (see
Sec.~1.1 and also~[3]).

       Note that the formulas $e_\alpha(u)\mapsto\lambda_\alpha(u)e_\alpha(u)$ define an automorphism of the algebra $X_1^m(\Gamma,\mu;\lambda)$
for any nonzero functions  $\lambda_1(u),\dots,\lambda_m(u)$. The commutative subalgebra consisting of meromorphic
functions of the variables $y_1, \dots, y_m$ is an analog of the
Cartan subalgebra.

       In this paper, we construct similar homomorphisms for arbitrary $n$ and~$k$. In this case, it turns
out that the structure of the algebras into which the algebra
$Z_{n,k}(\Gamma,\eta)$ is mapped depends on the expansion of the
number $\frac nk$ into a continued fraction. Namely, let $\frac
nk=n_1-\frac1{n_2-\ldots-\frac1{n_p}}$, where $n_\alpha\ge2$ and
$p$ is the length of the continued fraction. There is a family of algebras $X_p^{m_1,\dots,m_p}(\Gamma,\mu;\lambda_1,\dots,\lambda_p)$,
where $m_1, \dots, m_p\in   \N$ and
$\mu,\lambda_1,\dots,\lambda_p\in\C$. The algebra
$X_p^{m_1,\dots,m_p}(\Gamma,\mu;\lambda_1,\dots,\lambda_p)$
 is generated
by the commutative subalgebra consisting of meromorphic functions of
the variables $\bigl\{y_{\alpha,j};\,1\le j\le
p,\,\penalty50{}1\le\alpha\le m_j\bigr\}$ and by the generators
$\bigl\{e_{\alpha_1,\dots,\alpha_p}(u);\,1\le\alpha_j\le
m_j,\,u\in\C\bigr\}$. In this case, $y_{\alpha,j}$ are
dynamical variables, i.e., $e_{\alpha_1,\dots,\alpha_p}(u)y_{\beta,j}=
(y_{\beta,j}+\lambda_j)e_{\alpha_1,\dots,\alpha_p}(u)$ if $\beta\ne\alpha_j$, and
$e_{\alpha_1,\dots,\alpha_p}(u)y_{\alpha_j,j}=
(y_{\alpha_j,j}+\lambda_j-\mu)e_{\alpha_1,\dots,\alpha_p}(u)$. The
relations between $e_{\alpha_1,\dots,\alpha_p}(u)$ have the form
\begin{multline*}
e_{\alpha_1,\dots,\alpha_p}(u)e_{\beta_1,\dots,\beta_p}(v)\\
=\psi_1e_{\alpha_1,\dots,\alpha_p}(v)e_{\beta_1,\dots,\beta_p}(u)+
\psi_2e_{\beta_1,\alpha_2,\dots,\alpha_p}(v)
e_{\alpha_1,\beta_2,\dots,\beta_p}(u)+
\ldots\\
{}+\psi_pe_{\beta_1,\dots,\beta_{p-1},\alpha_p}(v)
e_{\alpha_1,\dots,\alpha_{p-1},\beta_p}(u)+
\psi_{p+1}e_{\beta_1,\dots,\beta_p}(v)e_{\alpha_1,\dots,\alpha_p}(u).
\end{multline*}
Here  $\psi_1,\dots,\psi_{p+1}$  are some functions of the variables $u,v,y_{\alpha_1,1},\dots,y_{\alpha_p,p},y_{\beta_1,1},\dots,y_{\beta_p,p}$. In this
relation,  $\alpha_j\ne\beta_j$ for all $1\le j\le p$. If   $\alpha_j=\beta_j$ for some $j$, then there are similar relations
(see Sec.~3). The formulas $e_{\alpha,j}(u)\mapsto
h_{\alpha,j}(u)e_{\alpha,j}(u)$ define an automorphism of the algebra
$X_p^{m_1,\dots,m_p}(\Gamma,\mu,\lambda_1,\dots,\lambda_p)$ for any nonzero functions $h_{\alpha,j}(u)$. For any $m_1, \dots , m_p\in   \N$, there is a
homomorphism $Z_{n,k}(\Gamma,\eta)\to
X_p^{m_1,\dots,m_p}(\Gamma,\mu;\lambda_1,
\dots,\lambda_p)$ defined by the formula
$$
x_i(u)\mapsto\sum_{\begin{subarray}{c} 1\le\alpha_1\le
m_1\\ \dots\\ 1\le\alpha_p\le
m_p\end{subarray}}w_i(y_{\alpha_1,1}+\nu_1u,\dots,y_{\alpha_p,p}+\nu_pu)
e_{\alpha_1,\dots,\alpha_p}(nu).
$$
Here  $\nu_1,\dots,\nu_p,\mu,\lambda_1,\dots,\lambda_p$ are expressed
in terms of $\eta$  and $n_1, \dots , n_p$, where $\frac
nk=n_1-\frac1{n_2-\ldots-\frac1{n_p}}$ 
(see Sec.~4), and $w_i(u_1,\dots,u_p)\in\Theta_{n/k}(\Gamma)$ are
theta functions of $p$ variables (see Sec.~1.2).

       There is a dual construction in which an algebra with similar exchange relations is mapped
into the algebra~$Z_{n,k}(\Gamma,\eta)$. In this case, the structure of the former algebra depends on the
expansion of the number $\frac n{n-k}$ into a continued
fraction. Namely, let  $\frac
n{n-k}=n_1'-\frac1{n_2'-\ldots-\frac1{n_{p'}'}}$, 
where $n_\alpha'\ge2$ and $p'$  is the length of the continued fraction. There is a family of algebras
$Y_{p'}(\Gamma,\mu;\mu_1,\dots,\mu_{p'})$, where $\mu,\mu_1,\dots,\mu_{p'}\in\C$. The
generators of the algebra 
are $\bigl\{e(u,u_1,\dots,u_{p'});\,u,u_1,\dots,u_{p'}\in\C\bigr\}$, and its relations have the form
\begin{multline*}
e(u,u_1,\dots,u_{p'})e(v,v_1+\mu_1,\dots,v_{p'}+\mu_{p'})\\
=f_1
e(v,u_1,\dots,u_{p'})e(u,v_1+\mu_1,\dots,v_{p'}+\mu_{p'})\\
{}+f_2
e(v,v_1,u_2,\dots,u_{p'})e(u,u_1+\mu_1,v_2+\mu_2,\dots,v_{p'}+\mu_{p'})+\ldots\\
{}+f_{p'+1}
e(v,v_1,\dots,v_{p'})e(u,u_1+\mu_1,\dots,u_{p'}+\mu_{p'}),
\end{multline*}
where $f_1,\dots,f_{p'+1}$ are some functions of
$u,u_1,\dots,u_{p'};v,v_1,\dots,v_{p'}$  (see Sec.~5).

       There is a homomorphism
$Y_{p'}(\Gamma,\mu,\mu_1,\dots,\mu_{p'})\to
Z_{n,k}(\Gamma,\eta)$ defined by the formula
$$
e(nu,u_1,\dots,u_{p'})\mapsto
\sum_{\alpha\in\Z/n\Z}x_{1-\alpha}(u)w_\alpha(u_1+\gamma_1u,
\dots,u_{p'}+\gamma_{p'}u).
$$
Here $\mu,\mu_1,\dots,\mu_{p'},\gamma_1,\dots,\gamma_{p'}$  can be written in terms of  $\eta,n_1',\dots,n_{p'}'$  (see Sec.~6). The expressions
$w_\alpha(u_1,\dots,u_{p'})\in\Theta_{n/n-k}(\Gamma)$ are theta
functions of $p'$  variables (see Sec. 1.2).

       We now describe the content of the paper.

   In Sec.~1, we collect the main facts needed below and related to theta functions. In Sec.~1.1,
the space of theta functions of one variable is considered; we denote this space by  $\Theta_{n,c}(\Gamma)$. Geometrically, these are holomorphic sections of a linear bundle on an elliptic curve $\C/\Gamma$  with Chern
class $n$ where $c\in\C$ is a continuous parameter of the bundle. For more detail, see [6], although we
use somewhat different notation.

   In Sec.~1.2, the space $\Theta_{n/k}(\Gamma)$ of theta functions of $p$ variables is introduced, where $p$ is the
length of the continued fraction $\frac
nk=n_1-\frac1{n_2-\ldots-\frac1{n_p}}$, $n_\alpha\ge2$.
These are sections of some special 
linear bundle on~$(\C/\Gamma)^p$. It can also be shown (but we do not
need this) that the space $\Theta_{n/k}(\Gamma)$ is
isomorphic to the space of sections of a holomorphic indecomposable
bundle on an elliptic curve~$\C/\Gamma$  
of degree $n$ and rank~$k$. This bundle is known to be unique to
within a shift on the elliptic curve.

   In Sec.~2, the Belavin $R$-matrix and the related Zamolodchikov
algebra $Z_{n,k}(\Gamma,\eta)$ are described.

   In Sec.~3, dynamical algebras
$X_p^{m_1,\dots,m_p}(\Gamma,\mu;\lambda_1,\dots,\lambda_p)$ with
exchange relations are constructed.

   In Sec.~4, homomorphisms of the Zamolodchikov algebra
$Z_{n,k}(\Gamma,\eta)$ into the dynamical exchange
algebras $X_p^{m_1,\dots,m_p}(\Gamma,\mu;\lambda_1,\dots,\lambda_p)$
are constructed.

   In Sec.~5, polyspectral exchange algebras
$Y_{p'}(\Gamma,\mu;\mu_1,\dots,\mu_{p'})$ are constructed.

   In Sec.~6, a homomorphism from the algebra
$Y_{p'}(\Gamma,\mu;\mu_1,\dots,\mu_{p'})$ into the Zamolodchikov
algebra $Z_{n,k}(\Gamma,\eta)$ is constructed.

   Here, it is found that, for a fixed $u$, the space of the
generators $\bigl\{x_i(u);\, i\in\Z/n\Z,\,u\in\C\bigr\}$ of
the algebra $Z_{n,k}(\Gamma,\eta)$ is naturally isomorphic to the
space $\Theta_{n/k}(\Gamma)$ of theta functions and is dual
to the space~$\Theta_{n/n-k}(\Gamma)$. The related duality between
these spaces is described in Sec.~1.3.

   The main results of this paper are contained in Propositions 6 and~7. In Proposition 6, homomorphisms from the Zamolodchikov algebra
$Z_{n,k}(\Gamma,\eta)$ into algebras with exchange relations and
dynamical variables are constructed. In Proposition 7, a homomorphism from a similar algebra with
exchange relations, but without dynamical variables, into the
Zamolodchikov algebra $Z_{n,k}(\Gamma,\eta)$ is
constructed. These two propositions follow from identity (2) proved in Sec.~1.2.

\section{Theta Functions Associated\\ with the Degree of an Elliptic Curve}

\subsection{Theta functions of one variable} 
Let $\Gamma\subset\C$ be the integer lattice generated by 1 and
$\tau\in\C$, where $\Im\tau>0$. Let $n\in\N$ and $c\in\C$. We denote
by  $\Theta_{n,c}(\Gamma)$ the space of entire functions
of one variable that satisfy the relations
$$
f(z+1)=f(z),\quad f(z+\tau)=(-1)^ne^{-2\pi i(nz-c)}f(z).
$$
It follows that
$$
f(z+\alpha+\beta\tau)=(-1)^{n\beta} e^{-2\pi
i\beta\left(nz-c+\frac{(\beta-1)n}2\tau\right)}f(z)
$$
for each element   $\alpha+\beta\tau\in\Gamma$       (i.e., for  $\alpha,\beta\in\Z$). As is known [6], $\dim\Theta_{n,c}(\Gamma)=n$, and every
function $f\in\Theta_{n,c}(\Gamma)$ has exactly $n$ zeros modulo~$\Gamma$   (counting multiplicity), whose sum is equal
to $c$ modulo~$\Gamma$. Let  $\theta(z)=\sum_{\alpha\in\Z}(-1)^\alpha
e^{2\pi i\left(\alpha
z+\frac{\alpha(\alpha-1)}2\tau\right)}$. It~is clear that
$\theta(z)\in\Theta_{1,0}(\Gamma)$. What has
been said implies that  $\theta(0)=0$ , and this is the only zero modulo~$\Gamma$. It~can readily be verified that
$\theta(-z)=-e^{-2\pi iz}\theta(z)$. Moreover, it is known that
$\theta(z)$ can be expanded as an infinite product [6],
$$
\theta(z)=\prod_{\alpha\ge1}(1-e^{2\pi
i\alpha\eta})\cdot(1-e^{2\pi
iz})\cdot\prod_{\alpha\ge1}(1-e^{2\pi
i(z+\alpha\eta)})(1-e^{2\pi i(\alpha\eta-z)}).
$$
We define linear operators $T_{\frac1n}$ and $T_{\frac1n\tau}$ acting on the space of functions of one variable,
$$
T_{\frac1n}f(z)=f\left(z+\frac1n\right),\quad
T_{\frac1n\tau}f(z)=e^{2\pi
i\left(z+\frac1{2n}-\frac{n-1}{2n}\tau\right)}f\left(z+\frac1n\tau\right).
$$
It can readily be seen that $T_{\frac1n}T_{\frac1n\tau}=e^{\frac{2\pi
i}n}T_{\frac1n\tau}T_{\frac1n}$. A simple calculation shows that the space
$\Theta_{n,\frac{n-1}2}(\Gamma)$ is invariant under $T_{\frac1n}$
and~$T_{\frac1n\tau}$. The restriction of $T_{\frac1n}$ and~$T_{\frac1n\tau}$ to  $\theta_{n,\frac{n-1}2}(\Gamma)$
satisfies the additional relations  $T_{\frac1n}^n=T_{\frac1n\tau}^n=1$. Let $G_n$ be the group with the generators $a, b$,
and $\epsilon$  and the relations $ab=\epsilon ba$, $a\epsilon=\epsilon
a$, $b\epsilon=\epsilon
b$, and $a^n=b^n=\epsilon^n=e$.  The group $G_n$ is
the central extension of the group $(\Z/n\Z)^2$, namely, the element
$\epsilon$   generates the normal subgroup
$C_n\cong\Z/n\Z$, and $G_n/C_n\cong(\Z/n\Z)^2$. The formulas $a\mapsto
T_{\frac1n}$ and $b\mapsto T_{\frac1n\tau}$ define an irreducible
representation of the group $G_n$ in the space   $\theta_{n,\frac{n-1}2}(\Gamma)$~[6]. We select a basis $\bigl\{\theta_\alpha,\,\alpha\in\Z/n\Z\bigr\}$
of~$\Theta_{n,\frac{n-1}2}(\Gamma)$  on which the above operators act
in the following way: $T_{\frac1n}\theta_\alpha=e^{2\pi i\frac\alpha n}\theta_\alpha$  and
$T_{\frac1n\tau}\theta_\alpha=\theta_{\alpha+1}$. It is clear that this can be done uniquely to within multiplication by a common
constant. Namely, we set
\begin{multline*}
\theta_\alpha(z)=\theta\left(z+\frac{\alpha}{n}\tau\right)
\theta\left(z+\frac{\alpha}{n}\tau+\frac1n\right)\dots
\theta\left(z+\frac{\alpha}{n}\tau+\frac{n-1}n\right)\\
{}\times
e^{2\pi i\left(\alpha
z+\frac{\alpha(\alpha-n)}{2n}\tau+\frac\alpha{2n}\right)}.
\end{multline*}
It is easy to show that
 $\theta_\alpha(z)\in\Theta_{n,\frac{n-1}2}(\Gamma)$,
$\theta_{\alpha+n}(z)=\theta_\alpha(z)$, and
\begin{align*}
\theta_\alpha\left(z+\frac1n\right)&=e^{2\pi i\frac\alpha
n}\theta_\alpha(z),\\
\theta_\alpha\left(z+\frac1n\tau\right)&=
e^{-2\pi i\left(z+\frac1{2n}-\frac{n-1}{2n}\tau\right)}\theta_{\alpha+1}(z).
\end{align*}

      Clearly, the functions
 $\left\{\theta_\alpha\left(z-\frac1nc-\frac{n-1}{2n}\right),
\alpha\in\Z/n\Z\right\}$ form a basis in~$\Theta_{n,c}(\Gamma)$.

      We shall need the identity
\begin{equation}
\theta(nz)=\frac{n\theta_0(z)\dots\theta_{n-1}(z)e^{-2\pi i\frac{n(n-1)}2z}}
{\theta_1(0)\dots\theta_{n-1}(0)\theta\left(\frac1n\right)\dots
\theta\left(\frac{n-1}n\right)}.
\end{equation}

\begin{proofd} 
of the identity. It can readily be verified that both sides of the identity belong to the space
$\Theta_{n^2,\frac{n(n-1)}2\tau}(\Gamma)$, and each of them has $n^2$
zeros modulo~$\Gamma$   at the points $\left\{\frac\alpha
n+\frac\beta n\tau;\alpha,\beta\in\Z\right\}$.
Consequently, they differ by a constant factor. The constant can readily be calculated by dividing
the two sides of the identity by~$\theta(z)$ and finding the limit as
$z\to 0$.
\end{proofd}

\begin{remark} 
The modular parameter $\tau$  is assumed to be fixed throughout the paper. It will
always be clear from the context what is the value of $n$ for~$\theta_\alpha(z)$. When necessary, we shall use the
notation~$\theta_\alpha^n(z)$  if the formula under consideration
involves several types of theta functions.
\end{remark}

\subsection{Theta functions of several variables} 

Let $n$ and $k$ be coprime positive integers such
that $1\le
k<n$. We expand $\frac nk$ into a continued fraction of the following
form: $\frac
nk=n_1-\frac1{n_2-\frac1{n_3-\ldots-\frac1{n_p}}}$,
where all $n_\alpha$  are greater than or equal to~2. It is clear that such an expansion exists and is
unique. Denote by $d(m_1, \dots , m_q)$ the determinant of the $q
\times q$ matrix~$(m_{\alpha\beta})$, where
$m_{\alpha\alpha}=m_\alpha$,
$m_{\alpha,\alpha+1}=m_{\alpha+1,\alpha}=-1$, and
$m_{\alpha,\beta}=0$ for $|\alpha-\beta|>1$. We set $d(\varnothing)=1$ for $q = 0$. It follows from
the elementary theory of continued fractions that $n = d(n_1, \dots ,
n_p)$ and $k = d(n_2, \dots , n_p)$.

      We again consider the integer lattice   $\Gamma\subset\C$
generated by 1 and~$\tau$, where $\Im\tau>0$.

      Let  $\Theta_{n/k}(\Gamma)$ be the space of entire functions $f$ of $p$ variables such that
\begin{align*}
f(z_1,\dots,z_\alpha+1,\dots,z_p)&=f(z_1,\dots,z_p),\\
f(z_1,\dots,z_\alpha+\tau,\dots,z_p)&=(-1)^{n_\alpha}e^{-2\pi
i(n_\alpha
z_\alpha-z_{\alpha-1}-z_{\alpha+1}-(\delta_{1,\alpha}-1)\tau)}
f(z_1,\dots,z_p).
\end{align*}
Here $1\le\alpha\le p$, $z_0=z_{p+1}=0$, and $\delta_{1,\alpha}$
is the Kronecker delta. Thus, the functions $f\in\Theta_{n/k}(\Gamma)$
are periodic with period 1 and quasiperiodic with period~$\tau$   in each of the variables. The periodicity
implies that each of the functions in $\Theta_{n/k}(\Gamma)$ can be expanded into a Fourier series of the form
$f(z_1,\dots,z_p)=\sum_{\alpha_1,\dots,\alpha_p\in\Z}
a_{\alpha_1\dots\alpha_p}e^{2\pi i(\alpha_1z_1+\ldots+\alpha_pz_p)}$. The quasiperiodicity gives a linear system of
equations for the coefficients,
$$
a_{\alpha_1,\dots,\alpha_{\nu-1}-1,\alpha_\nu+n_\nu,
\alpha_{\nu+1}-1,\dots,\alpha_p}=
(-1)^{n_\alpha}e^{2\pi i(\alpha_\nu+\delta_{1,\alpha}-1)\tau}
a_{\alpha_1\dots\alpha_p}.
$$

      It can readily be seen that this system has $n = d(n_1, \dots ,
n_p)$ linearly independent solutions
each defining a function in~$\Theta_{n/k}(\Gamma)$ (for $\Im\tau>0$,
$n_1,\dots,n_p\ge2$).

      For $k = 1$, we have the space
$\Theta_n(\Gamma)=\Theta_{n,0}(\Gamma)$ of functions of one variable in Sec.~1.1 with
the basis
$\left\{w_\alpha(z)=\theta_\alpha\left(z+\frac{n-1}2\right),\
\alpha\in\Z/n\Z\right\}$. A similar basis can be constructed in the space
$\Theta_{n/k}(\Gamma)$ for an arbitrary~$k$. We define operators
$T_{\frac1n}$ and~$T_{\frac1n\tau}$ in the space of functions of $p$
variables as follows:
\begin{align*}
T_{\frac1n}f(z_1,\dots,z_p)&=f(z_1+r_1,\dots,z_p+r_p),\\
T_{\frac1n\tau}f(z_1,\dots,z_p)&=e^{2\pi i(z_1+\phi)}
f(z_1+r_1\tau,\dots,z_p+r_p\tau).
\end{align*}
Here $r_\alpha=\frac{d(n_{\alpha+1},\dots,n_p)}{d(n_1,\dots,n_p)}$ and  $\phi\in\C$ is a constant.

    Clearly, $T_{\frac1n}T_{\frac1n\tau}=e^{2\pi i\frac
kn}T_{\frac1n\tau}T_{\frac1n}$. As in the case of theta functions of one variable, the space
 $\Theta_{n/k}(\Gamma)$ is invariant under the operators $T_{\frac1n}$
and $T_{\frac1n\tau}$, and the restriction of these operators
to~$\Theta_{n/k}(\Gamma)$ satisfies the relations $T_{\frac1n}^n=1$ and $T_{\frac1n\tau}^n=\mu$, where $\mu\in\C$. We choose a $\phi$  such that
$\mu=1$. This can obviously be done to within multiplication of
$T_{\frac1n\tau}$ by an $n$th root of unity.

\begin{statement} 
There is a basis $\bigl\{w_\alpha(z_1,\dots,z_p);\,\alpha\in\Z/n\Z\bigr\}$ in  $\Theta_{n/k}(\Gamma)$ such that
$$
T_{\frac1n}w_\alpha=e^{2\pi i\frac kn\alpha}w_\alpha,\quad
T_{\frac1n\tau}w_\alpha=w_{\alpha+1}.
$$
It is defined uniquely to within a constant factor.
\end{statement}

\begin{proof}
 Let $f\in\Theta_{n/k}(\Gamma)$ be an eigenvector of $T_{\frac1n}$ with eigenvalue~$\lambda$. Since $T_{\frac1n}^n=1$ on~$\Theta_{n/k}(\Gamma)$,
we have  $\lambda^n=1$. Moreover, $T_{\frac1n}T_{\frac1n\tau}f=e^{2\pi
i\frac
kn}T_{\frac1n\tau}T_{\frac1n}f=e^{2\pi i\frac kn}\lambda
T_{\frac1n\tau}f$; hence, $T_{\frac1n\tau}f$ is again
an eigenvector with the eigenvalue $e^{2\pi i\frac
kn}\lambda$. Since $n$ and $k$ are coprime, the factor $e^{2\pi i\frac
kn}$ is a
primitive $n$th root of unity. Hence, $\bigl\{T_{\frac1n\tau}^\alpha
f;\,\alpha=0,1,\dots,n-1\bigr\}$
 are eigenvectors of~$T_{\frac1n}$ with
distinct eigenvalues, and any of the $n$th roots of unity is an
eigenvalue of some $T_{\frac1n\tau}^\alpha f$. Let $w_0$
satisfy the condition $T_{\frac1n}w_0=w_0$. We set
$w_\alpha=T_{\frac1n\tau}^\alpha w_0$. It is clear that
$T_{\frac1n}w_\alpha=e^{2\pi i\frac
kn\alpha}w_\alpha$ and
$T_{\frac1n\tau}w_\alpha=w_{\alpha+1}$. We also have
$w_{\alpha+n}=w_\alpha$ since $T_{\frac1n\tau}^n=1$
on~$\Theta_{n/k}(\Gamma)$.
\end{proof}

    Note that, as in the case of a theta function of one variable, the group $G_n$ acts irreducibly on
the space $\Theta_{n/k}(\Gamma)\colon a\mapsto T_{\frac1n}$, $b\mapsto
T_{\frac1n\tau}$, $\epsilon\mapsto\text{(multiplication by $e^{2\pi
i\frac kn}$})$.

\begin{remark}
 Let $L$ be the group of linear automorphisms acting on the space of functions of $p$
variables by the formula
$$
gf(z_1,\dots,z_p)=e^{2\pi
i(\phi_1z_1+\ldots+\phi_pz_p+\lambda)}f(z_1+\psi_1,\dots,z_p+\psi_p)
$$
for $g\in L$. Obviously, $L$ is a $(2p + 1)$-dimensional Lie
group. Let $L'\subset L$ be the subgroup of
transformations preserving the space~$\Theta_{n/k}(\Gamma)$, i.e.,
$L'=\bigl\{g\in
L;\,g(\Theta_{n/k}(\Gamma))=\Theta_{n/k}(\Gamma)\bigr\}$, and
let $L''\subset L'$  consist of the elements leaving
$\Theta_{n/k}(\Gamma)$ fixed, i.e., $L''=\bigl\{g\in
L';\,gf=f\text{ for all }f\in\Theta_{n/k}(\Gamma)\bigr\}$. It can be
shown that the quotient group $L'/L''=\wt G_n$  
 is generated by the elements
$T_{\frac1n}$ and $T_{\frac1n\tau}$ and by the multiplications by
constants.
\end{remark}

    We shall use the notation $w_\alpha^{n/k}(z_1,\dots,z_p)$ if it is unclear from the context what theta functions
are meant.

    We shall need an identity relating the theta functions in  $\Theta_{1,0}(\Gamma),\Theta_{n,\frac{n-1}2}(\Gamma)$, and~$\Theta_{n/k}(\Gamma)$,
\begin{align}
&\frac{\theta(y_1-z_1+nv-nu)}{\theta(nv-nu)\theta(y_1-z_1)}\nonumber\\
&{}\times w_\alpha(y_1+m_1u_,\dots,y_p+m_pu)
w_\beta(z_1+m_1v+l_1,\dots,z_p+m_p,v+l_p)\nonumber\\
&+\sum_{1\le t\le p-1}
\frac{\theta(z_t-y_t+y_{t+1}-z_{t+1})}
{\theta(z_t-y_t)\theta(y_{t+1}-z_{t+1})}\nonumber\\
&{}\times w_\alpha(z_1+m_1u,\dots,z_t+m_tu,y_{t+1}+m_{t+1}u,\dots,y_p+m_pu)
\nonumber\\ &{}\times
w_\beta(y_1+m_1v+l_1,\dots,y_t+m_tv+l_t,z_{t+1}+m_{t+1}v+l_{t+1},
\dots,z_p+m_pv+l_p)\nonumber\\
&{}+\frac{\theta(z_p-y_p+n\eta)}{\theta(z_p-y_p)\theta(n\eta)}
w_\alpha(z_1+m_1u_1,\dots,z_p+m_pu)
\nonumber\\ &{}\times w_\beta(y_1+m_1v+l_1,\dots,y_p+m_pv+l_p)\nonumber\\
&{}=\frac1n\theta\left(\frac1n\right)\dots\theta\left(\frac{n-1}n\right)
\nonumber\\ &{}\times
\sum_{r\in\Z/n\Z}\frac{\theta_{\beta-\alpha+ r(k-1)}(v-u+\eta)}
{\theta_{ r k}(\eta)\theta_{\beta-\alpha- r}(v-u)}
w_{\beta- r}(y_1+m_1v,\dots,y_p+m_pv)
\nonumber\\ &{}\times
w_{\alpha+ r}(z_1+m_1u+l_1,\dots,z_p+m_pu+l_p).
\end{align}
Here $m_\alpha=d(n_{\alpha+1},\dots,n_p)\eta$  and $l_\alpha=d(n_1,\dots,n_{\alpha-1})\eta$.

\begin{proof}
 Denote by   
$\phi_{\alpha,\beta}(\eta,u,v,y_1,\dots,y_p,z_1,\dots,z_p)$
the difference between the right- and left-
hand sides in (2). A calculation shows that it satisfies the following relations:
\begin{equation}
\begin{aligned}
\phi_{\alpha,\beta}(\eta,\dots,y_\alpha+1,\dots,z_p)&=
\phi_{\alpha,\beta}(\eta,\dots,z_p),\\
\phi_{\alpha,\beta}(\eta,\dots,y_\alpha+\tau,\dots,z_p)&=
-e^{-2\pi i(n_\alpha y_\alpha-y_{\alpha-1}-y_{\alpha+1}+\delta_{\alpha,1}v)}
\phi_{\alpha,\beta}(\eta,\dots,z_p),\\
\phi_{\alpha,\beta}(\eta,\dots,z_\alpha+1,\dots,z_p)&=
\phi_{\alpha,\beta}(\eta,\dots,z_p),\\
\phi_{\alpha,\beta}(\eta,\dots,z_\alpha+\tau,\dots,z_p)&=
-e^{-2\pi i(n_\alpha z_\alpha-z_{\alpha-1}-z_{\alpha+1}+
\delta_{\alpha,1}u+\delta_{\alpha,p}\eta)}
\phi_{\alpha,\beta}(\eta,\dots,z_p).
\end{aligned}
\end{equation}
Here $y_0=y_{p+1}=z_0=z_{p+1}=0$ and   $\delta_{\alpha,\beta}$  is the Kronecker delta. Moreover, calculations show
that there are no poles on the divisors $nv-nu\in\Gamma$,
$n\eta\in\Gamma$,
$y_1-z_1\in\Gamma$, \dots, and $y_p-z_p\in\Gamma$,
and hence the function   $\phi_{\alpha,\beta}$   is holomorphic everywhere on~$\C^{2p+3}$. However, it is clear that the
functions $\bigl\{w_\lambda(y_1+m_1v,\dots,y_p+m_pv)
w_\nu(z_1+m_1u+l_1,\dots,z_p+m_pu+l_p);\,
\lambda,\nu\in\Z/n\Z\bigr\}$ form
a basis in the space of holomorphic functions of $y_1,\dots,y_p,z_1,\dots,z_p$ that satisfy conditions (3).
Therefore, the function   $\phi_{\alpha,\beta}$  has the form
\begin{multline}
\varphi_{\alpha,\beta}(\eta,u,v,y_1,\dots,z_p)\\
{}=\sum_{\lambda,\nu\in\Z/n\Z}
\psi_{\lambda,\nu}(\eta,u,v)w_\lambda(y_1+m_1v,\dots,y_p+m_pv)
\\ {}\times
w_\nu(z_1+m_1u+l_1,\dots,z_p+m_pu+l_p).
\end{multline}

        Here the functions   $\psi_{\lambda,\nu}(\eta,u,v)$ are holomorphic and satisfy the relations
\begin{equation}
\begin{gathered}
\psi_{\lambda,\nu}(\eta+1,u,v)=\psi_{\lambda,\nu}(\eta,u+1,v)=
\psi_{\lambda,\nu}(\eta,u,v+1)=\psi_{\lambda,\nu}(\eta,u,v),\\
\psi_{\lambda,\nu}(\eta+\tau,u,v)=
e^{-2\pi in(v-u)}\psi_{\lambda,\nu}(\eta,u,v),\\
\psi_{\lambda,\nu}(\eta,u+\tau,v)=e^{2\pi in\eta}\psi_{\lambda,\nu}(\eta,u,v),\\
\psi_{\lambda,\nu}(\eta,u,v+\tau)=e^{-2\pi in\eta}\psi_{\lambda,\nu}(\eta,u,v).
\end{gathered}
\end{equation}

        These relations can be verified by a calculation, namely, the multipliers under the shifts by 1
and~$\tau$   in formulas (3) and (4) should be compared.

        However, every holomorphic function of the variables $\eta$, $u$ and $v$ that satisfies relations (5) is
equal to 0. Indeed, the periodicity implies the expansion into the Fourier series,
$$
\psi_{\lambda,\nu}(\eta,u,v)=\sum_{\alpha,\beta,\gamma\in\Z}
a_{\lambda,\nu,\alpha,\beta,\gamma}
e^{2\pi i(\alpha\eta+\beta u+\gamma v)}.
$$
Furthermore, it follows from the quasiperiodicity that the
coefficients $a_{\lambda,\nu,\alpha,\beta,\gamma}$  are equal to 0.
\end{proof}

\subsection{Duality between the spaces  $\Theta_{n/k}(\Gamma)$ and
$\Theta_{n/n-k}(\Gamma)$} 

Let us construct the canonical element 
$\Delta_{n,k}\in\Theta_{n/k}(\Gamma)\otimes\Theta_{n/n-k}(\Gamma)$
realizing the duality between these spaces (see (6)).

\begin{statement}
 Let
$$
\frac nk=n_1-\frac1{n_2-\ldots-\frac1{n_p}},\quad \frac
n{n-k}=n_1'-\frac1{n_2'-\ldots-\frac1{n_{p'}'}}
$$ 
be expansions into continued fractions, where $n_\alpha\ge2$ and
$n_\beta'\ge2$  for $1\le\alpha\le p$ and $1\le\beta\le p'$.
Here $p$ and $p'$  are the lengths of the continued fractions. In this case, $p'=n_1+\ldots+n_p-2p+1$
and $n_1'+\ldots+n_{p'}'=2(n_1+\ldots+n_p)-3p+1$. Moreover,
$n_1'+\ldots+n_\alpha'=2\alpha+\beta$   for
$n_1+\ldots+n_\beta-2\beta+1\le\alpha\le
n_1+\ldots+n_{\beta+1}-2\beta-2$. In other words, the Young diagrams for the partitions
$(n_1-1,n_1+n_2-3,\dots,n_1+\ldots+n_\alpha-2\alpha+1,\dots)$ and $(n_1'-1,n_1'+n_2'-3,\dots,n_1'+\dots+n_\beta'-2\beta+1,\dots)$
are mutually dual.
\end{statement}

\begin{remark} 
We have $p'  = n - 1$ and $n_1'=\ldots=n_{n-1}'=2$ for $k = 1$, $p = 1$, and $n_1 = n$.
For $p > 1$, if $n_2, \dots , n_{p-1}\ge   3$, then the sequence
$(n_1',\dots,n_p')$ has the form
$(2^{(n_1-2)},3,2^{(n_2-3)},3,\dots,3,2^{(n_{p-1}-3)},3,2^{(n_p-2)})$.
Here $2^{(t)}$, $t\ge0$, denotes the sequence consisting of $t$ twos. This
formula is also true without the condition $n_2,\dots,n_p\ge3$ under the following convention: the
sequence $(m_1,2^{(-1)},m_2)$ has unit length and is equal to $(m_1+m_2-2)$. This rule is consecutively
applied to all $n_\alpha=2$ for $2\le\alpha\le p-1$.
\end{remark}

\begin{proof}
 The proof is carried out by induction on $\min(p,p')$. For $p = 1$, we must prove that
$\frac n{n-1}=2-\frac1{2-\ldots-\frac12}$  is a continued fraction of length $n - 1$. For example, let $n_1>2$ for $p,p'>1$. We
have     $\frac
k{d(n_3,\dots,n_p)}=n_2-\frac1{n_3-\ldots-\frac1{n_p}}$. By assumption,
$$
\frac k{k-d(n_3,\dots,n_p)}=
n_{n_1-1}'-1-\frac1{n_{n_1}'-\frac1{n_{n_1+1}'-\ldots-\frac1{n_{p'}'}}}.
$$
where the sequence $(n_1',\dots,n_{n_1-2}')$ is~$(2^{(n_1-2)})$.
We can further see that
$n_1'-\frac1{n_2'-\ldots-\frac1{n_{p'}'}}=\frac n{n-k}$ by
using the relations 
$d(n_1,\dots,n_p)=n$,
$d(n_2,\dots,n_p)=k$,
$$
\frac{d(n_1,\dots,n_p)}{d(n_2,\dots,n_p)}=
n_1-\frac1{n_2-\ldots-\frac1{n_p}},
$$
and $d(n_1,\dots,n_p)=n_1d(n_2,\dots,n_p)-d(n_3,\dots,n_p)$.
\end{proof}

\begin{statement} 
Let the function  $\Delta_{n.k}(z_1,\dots,z_p;z_1',\dots,z_{p'}')$
 of $p+p'$  variables $z_1,\dots,z_p,z_1',\dots,z_{p'}'$
  be defined by the formula
\begin{multline*}
\Delta_{n,k}(z_1,\dots,z_p,z_1',\dots,z_{p'}')\\
{}=e^{2\pi
iz_1'}\theta(z_1-z_1')\theta(z_p+z_{p'}')\cdot\prod_{1\le\alpha\le
p'-1}\theta(z_\alpha'-z_{\alpha+1}'+z_{n_1'+\ldots+n_\alpha'-2\alpha+1})
\\
{}\times\prod_{1\le\beta\le
p-1}\theta(z_\beta-z_{\beta+1}+z_{n_1+\ldots+n_\beta-2\beta+1}').
\end{multline*}
This function satisfies the following relations:
$$
\Delta_{n,k}(z_1,\dots,z_\alpha+1,\dots,z_{p'}')=
\Delta_{n,k}(z_1,\dots,z_\beta'+1,\dots,z_{p'}')=
\Delta_{n,k}(z_1,\dots,z_{p'}'),\\
$$
\begin{multline*}
\Delta_{n,k}(z_1,\dots,z_\alpha+\tau,\dots,z_{p'}')\\
=(-1)^{n_\alpha}e^{-2\pi i
(n_\alpha z_\alpha-z_{\alpha-1}-z_{\alpha+1}-(\delta_{\alpha,1}-1)\tau)}
\Delta_{n,k}(z_1,\dots,z_{p'}'),
\end{multline*}
\begin{multline*}
\Delta_{n,k}(z_1,\dots,z_\beta'+\tau,\dots,z_{p'}')\\
=(-1)^{n_\beta'}e^{-2\pi i
(n_\beta'z_\beta'-z_{\beta-1}'-z_{\beta+1}'-(\delta_{\beta,1}-1)\tau)}
\Delta_{n,k}(z_1,\dots,z_{p'}'),
\end{multline*}
where $z_0=z_{p+1}=z_0'=z_{p'+1}'=0$ and $\delta_{\alpha,1}$  is the
Kronecker delta.
\end{statement}

      The proof follows directly from the above description of the duality between the sequences
$(n_1,\dots,n_p)$ and $(n_1',\dots,n_{p'}')$.

\begin{statement}                                                           
\begin{equation}
\Delta_{n,k}(z_1,\dots,z_p;z_1',\dots,z_{p'}')=c_{n,k}
\sum_{\alpha\in\Z/n\Z}w_\alpha^{n/k}(z_1,\dots,z_p)
w_{1-\alpha}^{n/n-k}(z_1',\dots,z_{p'}'),
\end{equation}
where $c_{n,k}\in\C$ is some constant.
\end{statement}

\begin{proof} 
It follows from Proposition 3 that the function of the variables $z_1,\dots,z_p$ defined by
$\Delta_{n,k}$ belongs to the space~$\Theta_{n/k}(\Gamma)$. Similarly, the function of the variables $z_1',\dots,z_{p'}'$  defined by
$\Delta_{n,k}$  belongs to~$\Theta_{n/n-k}(\Gamma)$. Therefore,  
$$
\Delta_{n,k}(z_1,\dots,z_p;z_1',\dots,z_{p'}')=
\sum_{\alpha,\beta\in\Z/n\Z}\lambda_{\alpha,\beta}
w_\alpha^{n/k}(z_1,\dots,z_p)
w_\beta^{n/n-k}(z_1',\dots,z_{p'}').
$$
However, one can readily see that
$$
\Delta_{n,k}(z_1+ r_1,\dots,z_p+ r_p;z_1'+ r_1',\dots,z_{p'}'+ r_{p'}')=
e^{\frac{2\pi i}n}\Delta_{n,k}(z_1,\dots,z_{p'}'),
$$
where $ r_\alpha=\frac{d(n_1,\dots,n_{\alpha-1})}n$ and $
r_\beta'=\frac{d(n_1',\dots,n_{\beta-1}')}n$. It follows that
$\lambda_{\alpha,\beta}=0$ for  $\alpha+\beta\not\equiv1\mmod n$
(since $w_\alpha(z_1+ r_1,\dots,z_p+ r_p)=e^{2\pi i\frac\alpha
n}w_\alpha(z_1,\dots,z_p)$ and $w_\beta(z_1'+ r_1',\dots,z_{p'}'+
r_{p'}')=e^{2\pi i\frac\beta
n}w_\beta(z_1',\dots,z_{p'}')$). Hence,   $\lambda_{\alpha,\beta}=\lambda_\alpha\delta_{\alpha+\beta,1}$. Similarly,
\begin{multline*}
\Delta_{n,k}(z_1+ r_1\tau,\dots,z_p+ r_p\tau;z_1'+ r_1'\tau,
\dots,z_{p'}'+ r_{p'}'\tau)\\
{}=e^{2\pi i\left(\frac1n\tau-z_p-z_{p'}'\right)}
\Delta_{n,k}(z_1,\dots,z_{p'}').
\end{multline*}
This implies that $\lambda_\alpha=\lambda_{\alpha+1}$, i.e.,
$\lambda_\alpha$    does not depend on~$\alpha$.
\end{proof}

\section{Belavin $R$-Matrix}

\subsection{Main definitions} 

Let $V$ be an $n$-dimensional linear space, let $R(u, v)$ be a meromorphic
function of complex variables $u$ and $v$ with range in the set of
linear operators in $V\otimes V$ such that
$R(u, v)$ satisfies the condition $R(u, v)R(v, u) = 1$, and let
$R_{\alpha\beta}^{\gamma\delta}(u,v)$ be the matrix elements of the
function $R(u, v)$ with respect to some basis in~$V$.

       We denote by $A_R$ the associative algebra with the generators $\bigl\{x_\alpha(u);\,\alpha=1,\dots,n;\,u\in\C\bigr\}$ and
the defining relations                                        
$$
x_\alpha(u)x_\beta(v)=\sum_{\gamma,\delta}R_{\alpha\beta}^{\gamma\delta}(u,v)x_\delta(v)x_\gamma(u).
$$

\begin{definition}  
An operator-valued function $R(u, v)$ is called an \emph{$R$-matrix} if, for any $u, v$ and
$w$ in general position, the elements
$\bigl\{x_\alpha(u)x_\beta(v)x_\gamma(w);\,\alpha,\beta,\gamma=1,\dots,n\bigr\}$
 of the algebra $A_R$ are
linearly independent. In this case, $A_R$ is called the Zamolodchikov
algebra for the $R$-matrix $R(u, v)$.
\end{definition}

       As is known, for this condition to hold it is necessary and
sufficient that $R(u, v)$ satisfy the
Yang--Baxter equation, which can be written as
\begin{equation}
\sum_{\mu,\nu,t}R_{\alpha\beta}^{\nu\mu}(u,v)R_{\nu\gamma}^{l
t}(u,w)R_{\mu
t}^{\phi\psi}(v,w)=\sum_{\mu,\nu,t}R_{\beta\gamma}^{\nu\mu}(v,w)
R_{\alpha\mu}^{t\psi}(u,w)R_{t\nu}^{l\phi}(u,v).
\end{equation}
in terms of matrix elements. Note that this equation is preserved if
$R(u, v)$ is multiplied by an
arbitrary function~$\phi(u,v)$.

       Krichever [7] classified the solutions of Eq.~(7) for $n = 2$. No classification is known for greater
values of~$n$.

\subsection{Belavin $R$-matrix} 

Belavin [2] constructed a family of $R$-matrices. In the above notation, the related Zamolodchikov algebras are written as follows. As above, let $n$ and $k$ be coprime
positive integers such that $1\le k<n$. Let    $\Gamma\subset\C$ be the lattice generated by 1 and~$\tau$, $\Im\tau>0$,
and let $\eta\in\C$.

       We define the algebra $Z_{n,k}(\Gamma,\eta)$ by the generators $\bigl\{x_\alpha(u);\,\alpha\in\Z/n\Z,\,u\in\C\bigr\}$ and by the
relations
\begin{multline}
\frac{\theta_1(0)\dots\theta_{n-1}(0)\theta_0(v-u+\eta)\dots
\theta_{n-1}(v-u+\eta)}
{\theta_0(\eta)\dots\theta_{n-1}(\eta)\theta_0(v-u)\dots
\theta_{n-1}(v-u)}x_\alpha(u)x_\beta(v)\\
{}=\sum_{ r\in\Z/n\Z}
\frac{\theta_{\beta-\alpha+ r(k-1)}(v-u+\eta)}
{\theta_{k r}(\eta)\theta_{\beta-\alpha- r}(v-u)}x_{\beta- r}(v)
x_{\alpha+ r}(u).
\end{multline}
Here  $\bigl\{\theta_\alpha(u);\,\alpha\in\Z/n\Z\bigr\}$ is a basis in
the space  $\Theta_{n,\frac{n-1}2}(\Gamma)$.

\begin{statement}
$Z_{n,k}(\Gamma,\eta)$ is a Zamolodchikov algebra.
\end{statement}

\begin{proof} 
It is easy to verify the corresponding Yang--Baxter equation (7) directly. We substitute
$$
R_{\alpha\beta}^{\delta\gamma}(u,v)=\delta_{\alpha+\beta,\gamma+\delta}
\frac{\theta_{\beta-\alpha+(\beta-\gamma)(k-1)}(v-u+\eta)}
{\theta_{k(\beta-\gamma)}(\eta)\theta_{\beta-\delta}(v-u)}
$$
into (7) and compare the poles on the left- and right-hand sides. Let
$\phi(u,v,w,\eta)$ be the difference between the left- and right-hand
sides. It can readily be shown that  $\phi(u,v,w,\eta)$ is holomorphic
in $u$ and satisfies the conditions
$\phi(u+1,v,w,\eta)=\phi(u,v,w,\eta)$ and
$\phi(u+\tau,v,w,\eta)=e^{2\pi in(2\eta)}\phi(u,v,w,\eta)$, and this
implies that $\phi=0$.
\end{proof}

\begin{remark} 
Suppose that $n\eta\in\Gamma$, i.e.,  $\eta=\frac\mu n+\frac\nu n\tau$, where  $\mu,\nu\in\Z$. Then relations (8)
become
$$
x_\alpha(u)x_\beta(v)=e^{\frac{2\pi i}n(\beta-\alpha+k'\nu)\mu}
x_{\beta+k'\nu}(v)x_{\alpha-k'\nu}(u).
$$
where $k'=d(n_1,\dots,n_{p-1})$, and it one can check that $1\le k'<n$ and
$kk'\equiv1\mmod n$.
\end{remark}

\section{Dynamical Algebras with Exchange\\ Relations}

       Let $p\in\N$, let $\mu,\lambda_1,\dots,\lambda_p\in\C$, and let $m_1,\dots,m_p\in\N$. We define an associative algebra
$X_p^{m_1,\dots,m_p}(\Gamma,\mu;\lambda_1,\dots,\lambda_p)$ in the following way: it is generated by the commutative subalgebra
consisting of all meromorphic functions of the variables $\bigl\{y_{\alpha,j};\,1\le j\le
p,\,1\le\alpha\le m_j,\,\alpha,j\in\N\bigr\}$ and
by the generators
$\bigl\{e_{\alpha_1,\dots,\alpha_p}(u);\,\alpha_1,\dots,
\alpha_p\in\N,\,1\le\alpha_j\le m_j,\,u\in\C\bigr\}$. The defining relations are
\begin{equation}
\begin{aligned}
e_{\alpha_1,\dots,\alpha_p}(u)y_{\beta,j}&=
(y_{\beta,j}+\lambda_j)e_{\alpha_1,\dots,\alpha_p}(u),
\qquad\text{where $\beta\ne\alpha_j$},\\
e_{\alpha_1,\dots,\alpha_p}(u)y_{\alpha_j,j}&=
(y_{\alpha_j,j}+\lambda_j-\mu)e_{\alpha_1,\dots,\alpha_p}(u).
\end{aligned}
\end{equation}
This means that $y_{\beta,j}$ are dynamical variables. The remaining relations are quadratic with respect
to~$e_{\alpha_1,\dots,\alpha_p}(u)$. We first write the relations in ``general position.'' Suppose that  $\alpha_1\ne\beta_1$, \dots,
$\alpha_p\ne\beta_p$.
Then
\begin{multline*}
\frac{\theta(v-u+\mu)}{\theta(v-u)}e_{\alpha_1,\dots,\alpha_p}(u)
e_{\beta_1,\dots,\beta_p}(v)\\
=\frac{\theta(\mu)\theta(v-u+y_{\alpha_1,1}-y_{\beta_1,1})}
{\theta(v-u)\theta(y_{\alpha_1,1}-y_{\beta_1,1})}
e_{\alpha_1,\dots,\alpha_p}(v) e_{\beta_1,\dots,\beta_p}(u)\\
{}+
\sum_{1\le t<p}\frac{\theta(\mu)\theta(y_{\alpha_t,t}-y_{\beta_t,t}+
y_{\alpha_{t+1},t+1}-y_{\beta_{t+1},t+1})}
{\theta(y_{\alpha_t,t}-y_{\beta_t,t})\theta(y_{\alpha_{t+1},t+1}-
y_{\beta_{t+1},t+1})}\\
{}\times
e_{\beta_1,\dots,\beta_t,\alpha_{t+1},\dots,\alpha_p}(v)
e_{\alpha_1,\dots,\alpha_t,\beta_{t+1},\dots,\beta_p}(u)\\
{}+
\frac{\theta(y_{\alpha_p,p}-y_{\beta_p,p}+\mu)}
{\theta(y_{\alpha_p,p}-y_{\beta_p,p})}
e_{\beta_1,\dots,\beta_p}(v)e_{\alpha_1,\dots,\alpha_p}(u).
\end{multline*}
The relations in nongeneral position appear if
$\alpha_\nu=\beta_\nu$    for some $\alpha_\nu$   and
$\beta_\nu$. Let  $\alpha_1\ne\beta_1$, \dots,
$\alpha_{\nu-1}\ne\beta_{\nu-1}$, and   $\alpha_\nu=\beta_\nu$,
$1\le\nu\le p$. Then $e_{\alpha_1,\alpha_2,\dots,\alpha_p}(u)
e_{\alpha_1,\beta_2,\dots,\beta_p}(v)=
e_{\alpha_1,\alpha_2,\dots,\alpha_p}(v)
e_{\alpha_1,\beta_2,\dots,\beta_p}(u)$ for $\nu=1$ and
\begin{multline*}
\frac{\theta(v-u+\mu)}{\theta(v-u)}
e_{\alpha_1,\dots,\alpha_{\nu-1},\alpha_\nu,\dots,\alpha_p}(u)
e_{\beta_1,\dots,\beta_{\nu-1},\alpha_\nu,\dots,\beta_p}(v)\\
=\frac{\theta(\mu)\theta(v-u+y_{\alpha_1,1}-y_{\beta_1,1})}
{\theta(v-u)\theta(y_{\alpha_1,1}-y_{\beta_1,1})}
e_{\alpha_1,\dots,\alpha_{\nu-1},\alpha_\nu,\dots,\alpha_p}(v)
e_{\beta_1,\dots,\beta_{\nu-1},\alpha_\nu,\dots,\beta_p}(u)\\
{}+
\sum_{1\le t<\nu-1}
\frac{\theta(\mu)\theta(y_{\alpha_t,t}-y_{\beta_t,t}+
y_{\alpha_{t+1},t+1}-y_{\beta_{t+1},t+1})}
{\theta(y_{\alpha_t,t}-y_{\beta_t,t})
\theta(y_{\alpha_{t+1},t+1}-y_{\beta_{t+1},t+1})}\\
{}\times 
e_{\beta_1,\dots,\beta_t,\alpha_{t+1},\dots,\alpha_{\nu-1},
\alpha_\nu,\dots,\alpha_p}(v)
e_{\alpha_1,\dots,\alpha_t,\beta_{t+1},\dots,\beta_{\nu-1},
\alpha_\nu,\dots,\beta_p}(u)\\
{}+
\frac{\theta(y_{\alpha_{\nu-1},\nu-1}-y_{\beta_{\nu-1},\nu-1}+\mu)}
{\theta(y_{\alpha_{\nu-1},\nu-1}-y_{\beta_{\nu-1},\nu-1})}
e_{\beta_1,\dots,\beta_{\nu-1},\alpha_\nu,\dots,\alpha_p}(v)
e_{\alpha_1,\dots,\alpha_{\nu-1},\alpha_\nu,\dots,\beta_p}(u),
\end{multline*}
for $\nu>1$. Here the subscripts in
$e_{\alpha_1,\dots,\alpha_p}(v)e_{\beta_1,\dots,\beta_p}(u)$ on the right-hand side with indices from 1
to $(\nu-1)$ are permuted, whereas the others remain fixed.

       Finally, let  $\alpha_\nu=\beta_\nu$,
$\alpha_\lambda=\beta_\lambda$, and   $\alpha_i\ne\beta_i$ for
$\nu<i<\lambda$. Here $1\le\nu<p$ and $\nu<\lambda\le p+1$
(the case $\lambda=p+1$ means that $\alpha_i\ne\beta_i$ for all $i>\nu$). Therefore, the relations
\begin{multline*}
e_{\alpha_1,\dots,\alpha_\nu,\alpha_{\nu+1},\dots,
\alpha_{\lambda-1},\dots,\alpha_p}(v)
e_{\beta_1,\dots,\alpha_\nu,\beta_{\nu+1},\dots,
\beta_{\lambda-1},\alpha_\lambda,\dots,\beta_p}(u)\\
=\sum_{\nu+1\le t<\lambda-1}
\frac{\theta(\mu)\theta(y_{\alpha_t,t}-y_{\beta_t,t}+
y_{\alpha_{t+1},t+1}-y_{\beta_{t+1},t+1})}
{\theta(y_{\alpha_t,t}-y_{\beta_t,t})
\theta(y_{\alpha_{t+1},t+1}-y_{\beta_{t+1},t+1})}\\
{}\times
e_{\alpha_1,\dots,\alpha_\nu,\beta_{\nu+1},\dots,\beta_t,
\alpha_{t+1},\dots,\alpha_p}(v)
e_{\beta_1,\dots,\alpha_\nu,\alpha_{\nu+1},\dots,\alpha_t,
\beta_{t+1},\dots,\beta_p}(u)\\
{}+
\frac{\theta(y_{\alpha_{\lambda-1},\lambda-1}-
y_{\beta_{\lambda-1},\lambda-1}+\mu)}
{\theta(y_{\alpha_{\lambda-1},\lambda-1}-
y_{\beta_{\lambda-1},\lambda-1})}\\
{}\times
e_{\alpha_1,\dots,\alpha_\nu,\beta_{\nu+1},\dots,\beta_{\lambda-1},
\dots,\alpha_p}(v)
e_{\beta_1,\dots,\alpha_\nu,\alpha_{\nu+1},\dots,\alpha_{\lambda-1},
\dots,\beta_p}(u).
\end{multline*}
hold. Here the subscripts with indices from ($\nu+1$) to ($\lambda-1$)
are permuted, whereas the others
remain fixed.

\begin{remark} 
Let $p = 1$. Then $X_1^m(\Gamma,\mu;\lambda)$ is the Zamolodchikov
algebra for a dynamical $R$-matrix~[5]. This means that
$X_1^m(\Gamma,\mu;\lambda)$ is a plane deformation of the ring of polynomials in infinitely
many variables $\bigl\{e_\alpha(u);\,1\le\alpha\le m,\,\penalty50u\in\C\bigr\}$ over the field of meromorphic functions of the variables $y_{1,1},\dots,y_{m,1}$. Under the deformation, this field of functions becomes a field of quasiconstants
(see (9)).

      The structure of the algebra $X_p^{m_1,\dots,m_p}(\Gamma,\mu;\lambda_1,\dots,\lambda_p)$
for $p > 1$ is more complicated. Namely,
for $\mu=\lambda_1=\ldots=\lambda_p=0$, this is a commutative algebra over the field of meromorphic functions of
the variables $\bigl\{y_{\alpha,t};\,1\le\alpha\le m_t,\,1\le t\le
p\bigr\}$, with the generators
$\bigr\{e_{\alpha_1,\dots,\alpha_p}(u);1\le\alpha_t\le
m_t,u\in\C\bigr\}$ 
and the relations
\begin{multline*}
e_{\alpha_1,\dots,\alpha_{\nu-1},\alpha_\nu,\dots,\alpha_p}(u)
e_{\beta_1,\dots,\beta_{\nu-1},\alpha_\nu,\dots,\beta_p}(v)\\
{}=
e_{\beta_1,\dots,\beta_{\nu-1},\alpha_\nu,\dots,\alpha_p}(u)
e_{\alpha_1,\dots,\alpha_{\nu-1},\alpha_\nu,\dots,\beta_p}(v)
\end{multline*}
for  $\alpha_1\ne\beta_1$, \dots, $\alpha_{\nu-1}\ne\beta_{\nu-1}$, and
\begin{multline*}
e_{\alpha_1,\dots,\alpha_\nu,\alpha_{\nu+1},\dots,
\alpha_{\lambda-1},\dots,\alpha_p}(v)
e_{\beta_1,\dots,\alpha_\nu,\beta_{\nu+1},\dots,
\beta_{\lambda-1},\alpha_\lambda,\dots,\beta_p}(u)\\
{}=
e_{\alpha_1,\dots,\alpha_\nu,\beta_{\nu+1},\dots,
\beta_{\lambda-1},\dots,\alpha_p}(v)
e_{\beta_1,\dots,\alpha_\nu,\alpha_{\nu+1},\dots,
\alpha_{\lambda-1},\dots,\beta_p}(u),
\end{multline*}
for    $\alpha_\nu=\beta_\nu$, $\alpha_\lambda=\beta_\lambda$ and
      $\alpha_i\ne\beta_i$,  $\nu<i<\lambda$.

      One can see that these relations admit a uniformization:
      $e_{\alpha_1,\dots,\alpha_p}(u)=
e_{\alpha_1}^1(u)e_{\alpha_1,\alpha_2}^2\dots
e_{\alpha_{p-1},\alpha_p}^p$,
where $\{e_\alpha^1(u),e_{\alpha_{i-1},\alpha_i}^i\}$ are independent
      variables. The algebra
      $X_p^{m_1,\dots,m_p}(\Gamma,\mu;\lambda_1,\dots,\lambda_p)$ is a plane
deformation of this commutative algebra.
\end{remark}

\section{Homomorphisms of the Algebras $Z_{n,k}(\Gamma,\eta)$ into Dynamical Algebras
                                                              with Exchange Relations}

\begin{statement}
 For an arbitrary sequence $m_1,\dots,m_p\in\N$, there is a homomorphism
$$
\Phi\colon Z_{n,k}(\Gamma,\eta)\to
X_p^{m_1,\dots,m_p}(\Gamma,\mu;\lambda_1,\dots,\lambda_p),
$$
which is defined on the generators by the formula
\begin{equation}
\Phi(x_\alpha(u))=
\sum_{\begin{subarray}{c}1\le\alpha_1\le
m_1\\ \dots\\ 1\le\alpha_p\le m_p\end{subarray}}
w_\alpha(y_{\alpha_1,1}+\nu_1u,\dots,y_{\alpha_p,p}+\nu_pu)
e_{\alpha_1,\dots,\alpha_p}(nu).
\end{equation}
Here we use the notation
  $w_\alpha(y_1,\dots,y_p)\in\Theta_{n/k}(\Gamma)$,
$\nu_j=d(n_{j+1},\dots,n_p)\eta$,
$\mu=d(n_1,\dots,n_p)\eta=n\eta$,
$\lambda_j=d(n_1,\dots.n_{j-1})\eta$  for $1\le j\le p$, and $\frac
nk=n_1-\frac1{n_2-\ldots-\frac1{n_p}}$.
\end{statement}

\begin{remark} 
Formula (10) shows that, for a fixed $u\in\C$, the space of the generators of the algebra
$Z_{n,k}(\Gamma,\eta)$ is naturally isomorphic to the
space~$\Theta_{n/k}(\Gamma)$. Here any basis element $x_\alpha(u)$ corresponds
to $w_\alpha(y_1,\dots,y_p)\in\Theta_{n/k}(\Gamma)$.
\end{remark}

\begin{proof}[Proof of Proposition \emph6] 
We must show that the image of relations (8) under the homomorphism~$\Phi$   holds in the algebra $X_p^{m_1,\dots,m_p}(\Gamma,\mu;\lambda_1,\dots,\lambda_p)$.
Applying $\Phi$  to the difference between the
left- and right-hand sides in (8) and using the relations in the algebra $X_p^{m_1,\dots,m_p}(\Gamma,\mu;\lambda_1,\dots,\lambda_p)$,
we obtain an expression of the form
$$
\sum_{\begin{subarray}{c}1\le\alpha_1,\beta_1\le
m_1\\ \dots\\ 1\le\alpha_p,\beta_p\le m_p\end{subarray}}
\psi_{\alpha_1,\dots,\alpha_p,\beta_1,\dots,\beta_p}
(y_{\alpha_1,1},\dots,y_{\alpha_p,p},y_{\beta_1,1},\dots,y_{\beta_p,p})
e_{\alpha_1,\dots,\alpha_p}(nv)e_{\beta_1,\dots,\beta_p}(nu).
$$
We must prove that $\psi_{\alpha_1,\dots,\beta_p}
(y_{\alpha_1,1},\dots,y_{\beta_p,p})=0$ for all  $\alpha_1,\dots,\alpha_p$ and $\beta_1,\dots,\beta_p$. Let us verify
this formula for  $\alpha_1\ne\beta_1$, \dots, $\alpha_p\ne\beta_p$. The related calculation shows that
\begin{multline*}
\psi_{\alpha_1,\dots,\beta_p}(y_{\alpha_1,1},\dots,y_{\beta_p,p})\\
{}=
\frac
{\theta_1(0)\dots\theta_{n-1}(0)\theta_0(v-u+\eta)\dots\theta_{n-1}(v-u+\eta)}
{\theta_0(\eta)\dots\theta_{n-1}(\eta)\theta_0(v-u)\dots\theta_{n-1}(v-u)}
\cdot\frac{\theta(nv-nu+n\eta)}{\theta(nv-nu)}\\
\times
\biggl(w_\alpha(y_{\alpha_1,1}+\nu_1u,\dots,y_{\alpha_p,p}+\nu_pu)
w_\beta(y_{\beta_1,1}+\nu_1v+\lambda_1,\dots,y_{\beta_p,p}+\nu_pv+\lambda_p)\\
{}\times
\frac{\theta(n\eta)\theta(nv-nu+y_{\alpha_1,1}-y_{\beta_1,1})}
{\theta(nv-nu)\theta(y_{\alpha_1,1}-y_{\beta_1,1})}\\
{}+\sum_{1\le t<p}
\frac{\theta(n\eta)
\theta(y_{\beta_t,t}-y_{\alpha_t,t}+y_{\alpha_{t+1},t+1}-y_{\beta_{t+1},t+1})}
{\theta(y_{\beta_t,t}-y_{\alpha_t,t})
\theta(y_{\alpha_{t+1},t+1}-y_{\beta_{t+1},t+1})}\\
{}\times
w_\alpha(y_{\beta_1,1}+\nu_1u,\dots,y_{\beta_t,t}+\nu_tu,
y_{\alpha_{t+1},t+1}+\nu_{t+1}u,\dots,y_{\alpha_p,p}+\nu_pu)\\
{}\times
w_\beta(y_{\alpha_1,1}+\nu_1v+\lambda_1,\dots,
y_{\alpha_t,t}+\nu_tv+\lambda_t,
y_{\beta_{t+1},t+1}+\nu_{t+1}v+\lambda_{t+1},\dots,
y_{\beta_p,p}+\nu_pv+\lambda_p)\\
{}+
\frac{\theta(y_{\beta_p,p}-y_{\alpha_p,p}+n\eta)}
{\theta(y_{\beta_p,p}-y_{\alpha_p,p})}
\\
{}\times
w_\alpha(y_{\beta_1,1}+\nu_1u,\dots,y_{\beta_p,p}+\nu_pu)
w_\beta(y_{\alpha_1,1}+\nu_1v+\lambda_1,\dots,y_{\alpha_p,p}+\nu_pv+\lambda_p)\biggr)\\
{}-
\sum_{ r\in\Z/n\Z}
\frac{\theta_{\beta-\alpha+ r(k-1)}(v-u+\eta)}
{\theta_{k r}(\eta)\theta_{\beta-\alpha- r}(v-u)}
w_{\beta- r}(y_{\alpha_1,1}+\nu_1v,\dots,y_{\alpha_p,p}+\nu_pv)\\
{}\times
w_{\alpha+ r}(y_{\beta_1,1}+\nu_1u+\lambda_1,\dots,
y_{\beta_p,p}+\nu_pu+\lambda_p).
\end{multline*}
We replace  $\theta_0(v-u+\eta)\dots\theta_{n-1}(v-u+\eta)$,
$\theta_0(\eta)\dots\theta_{n-1}(\eta)$, and
$\theta_0(v-\nobreak u)\dots\theta_{n-1}(v-\nobreak u)$ in
this expression by using identity (1), after which the relation $\psi_{\alpha_1,\dots,\beta_p}=0$ readily follows from
identity (2).

    The case in which the relation $\alpha_j=\beta_j$  holds for some  $\alpha_j$ and  $\beta_j$ can be treated similarly.
\end{proof}

\section{Polyspectral Algebras with Exchange Relations}

    Let $p'\in\N$ and $\mu,\mu_1,\dots,\mu_{p'}\in\C$. Let
$Y_{p'}(\Gamma,\mu;\mu_1,\dots,\mu_{p'})$ be the associative algebra with the
generators
$\bigl\{e(u,u_1,\dots,u_{p'});\,u,u_1,\dots,u_{p'}\in\C\bigr\}$  and the defining relations
\begin{multline}
\frac{\theta(v-u+\mu)}{\theta(v-\mu)}
e(u,u_1,\dots,u_{p'})e(v,v_1+\mu_1,\dots,v_{p'}+\mu_{p'})\\
=\frac{\theta(\mu)\theta(v-u+u_1-v_1)}{\theta(v-u)\theta(u_1-v_1)}
e(v,u_1,\dots,u_{p'})e(u,v_1+\mu_1,\dots,v_{p'}+\mu_{p'})\\
{}+
\sum_{1\le t<p'}
\frac{\theta(\mu)\theta(v_t-u_t+u_{t+1}-v_{t+1})}
{\theta(v_t-u_t)\theta(u_{t+1}-v_{t+1})}
e(v,v_1,\dots,v_t,u_{t+1},\dots,u_{p'})\\
{}\times
e(u,u_1+\mu_1,\dots,u_t+\mu_t,v_{t+1}+\mu_{t+1},\dots,v_{p'}+\mu_{p'})\\
{}+
\frac{\theta(v_{p'}-u_{p'}+\mu)}{\theta(v_{p'}-u_{p'})}
e(v,v_1,\dots,v_{p'})e(u,u_1+\mu_1,\dots,u_{p'}+\mu_{p'}).
\end{multline}

\begin{remark} 
For $\mu=\mu_1=\ldots=\mu_{p'}=0$, the algebra
$Y_{p'}(\Gamma,0,\dots,0)$ is the ring of polynomials in
infinitely many variables $\bigl\{e(u,u_1,\dots,u_{p'});\,u,\dots,u_{p'}\in\nobreak\C\bigr\}$. It can be shown that $Y_{p'}(\Gamma,\mu;\mu_1,\dots,\mu_{p'})$
is a plane deformation of $Y_{p'}(\Gamma,0,\dots,0)$. This means that
$Y_{p'}(\Gamma,\mu;\mu_1,\dots,\mu_{p'})$ is the Zamolodchikov
algebra for an $R$-matrix in the space of functions of the variables $u_1,\dots,u_{p'}$, i.e., $u$ is the spectral
parameter for $e(u,u_1,\dots,u_{p'})$, and $u_1,\dots,u_{p'}$  index
the basis (that is, $u_1,\dots,u_{p'}$   are an analog
of $i$ for $x_i(u)$ in a finite-dimensional $R$-matrix). The corresponding Yang--Baxter equation can be
verified directly.
\end{remark}

\section{Homomorphism of a Polyspectral Algebra with Exchange Relations
                                                       into the
                                                       Algebra
                                                       $Z_{n,k}(\Gamma,\eta)$}

\begin{statement} 
There is a homomorphism
$$
\Psi\colon Y_{p'}(\Gamma,\mu;\mu_1,\dots,\mu_{p'})\to
Z_{n,k}(\Gamma,\eta),
$$
which is defined on the generators by the formula
\begin{equation}
\Psi(e(nu,u_1,\dots,u_{p'}))=\sum_{\alpha\in\Z/n\Z}x_{1-\alpha}(u)
w_\alpha(u_1+\gamma_1u,\dots,u_{p'}+\gamma_{p'}u).
\end{equation}
Here
\begin{gather*}
\frac n{n-k}=n_1'-\frac1{n_2'-\ldots-\frac1{n_{p'}'}},\quad
\mu=n\eta,\\ 
\mu_j=d(n_1',\dots,n_{j-1}')\eta,\quad
\gamma_j=-d(n_{j+1}',\dots,n_{p'}')\eta,\quad 1\le j\le p'.
\end{gather*}
\end{statement}

\begin{remark} 
Formula (12) shows that, if $u$ is fixed, then the space of generators
for the algebra $Z_{n,k}(\Gamma,\eta)$ is naturally dual to the
space~$\Theta_{n/n-k}(\Gamma)$ (see Sec. 1.3).
\end{remark}

\begin{proof}[Proof of Proposition \emph7] 
We must prove the validity of the image of relations (11) 
in the algebra~$Z_{n,k}(\Gamma,\eta)$ with respect to the homomorphism~$\Psi$. Applying $\Psi$  to the difference between the left-
and right-hand sides of (11) and using relations (8) in the
algebra~$Z_{n,k}(\Gamma,\eta)$, one can readily represent the
resulting expression in the form  $\sum_{\gamma,\delta\in\Z/n\Z}\psi_{\gamma,\delta}(u,v,
u_1,\dots,u_{p'},v_1,\dots,v_{p'})x_{1-\gamma}(v)x_{1-\delta}(u)$.
We must prove that   $\psi_{\gamma,\delta}(u,\dots,v_{p'})=0$. One can readily see that
\begin{multline*}
\psi_{\gamma,\delta}(u,\dots,v_{p'})\\
{}=\frac{\theta(nv-nu+n\eta)}{\theta(nv-nu)}
\sum_{ r\in\Z/n\Z}
w_{\delta- r}(u_1+\gamma_1u,\dots,u_{p'}+\gamma_{p'}u)\\
{}\times
w_{\gamma+ r}(v_1+\mu_1+\gamma_1v,\dots,v_{p'}+\mu_{p'}+\gamma_{p'}v)
\frac{\theta_{\delta-\gamma- r(k+1)}(v-u+\eta)}
{\theta_{k r}(\eta)\theta_{\delta-\gamma- r}(v-u)}\\
{}-\frac{\theta_1(0)\dots\theta_{n-1}(0)\theta_0(v-u+\eta)\dots
\theta_{n-1}(v-u+\eta)}
{\theta_0(\eta)\dots\theta_{n-1}(\eta)\theta_0(v-u)\dots
\theta_{n-1}(v-u)}\\
{}\times\biggl(\frac{\theta(n\eta)\theta(nv-nu+u_1-v_1)}
{\theta(nv-nu)\theta(u_1-v_1)}
w_\gamma(u_1+\gamma_1v,\dots,u_{p'}+\gamma_{p'}v)\\
{}\times
w_\delta(v_1+\mu_1+\gamma_1u,\dots,v_{p'}+\mu_{p'}+\gamma_{p'}u)\\
{}+
\sum_{1\le t<p'}
\frac{\theta(n\eta)\theta(v_t-u_t+u_{t+1}-v_{t+1})}
{\theta(v_t-u_t)\theta(u_{t+1}-v_{t+1})}\\
{}\times w_\gamma(v_1+\gamma_1v,\dots,v_t+\gamma_tv,
u_{t+1}+\gamma_{t+1}v,\dots,u_{p'}+\gamma_{p'}v)\\
{}\times
w_\delta(u_1+\gamma_1u+\mu_1,\dots,u_t+\mu_t+\gamma_tu,
v_{t+1}+\mu_{t+1}+\gamma_{t+1}u,\dots,v_p'+\mu_{p'}+\gamma_{p'}u)\\
{}+
\frac{\theta(v_{p'}-u_{p'}+n\eta)}{\theta(v_{p'}-u_{p'})}
w_\gamma(v_1+\gamma_1v,\dots,v_{p'}+\gamma_{p'}v)\\
{}\times
w_\delta(u_1+\mu_1+\gamma_1u,\dots,u_{p'}+\mu_{p'}+\gamma_{p'}u)\biggr).
\end{multline*}
We again replace  $\theta_0(v-u+\eta)\dots\theta_{n-1}(v-u+\eta)$,
$\theta_0(\eta)\dots\theta_{n-1}(\eta)$, and  $\theta_0(v-u)\dots\theta_{n-1}(v-u)$
according to (1). Moreover, we make the change of variables
$u_j\mapsto u_j-\gamma_ju-\gamma_jv$, $v_j\mapsto
v_j-\gamma_ju-\gamma_jv$.
This results in identity~(2), but, in this case, for the theta functions in the space~$\Theta_{n/n-k}(\Gamma)$.
\end{proof}

\begin{remarks} 
1. The composition of the homomorphisms $\Psi$  and $\Phi$  gives the homomorphism
$$
\Phi\circ\Psi\colon Y_{p'}(\Gamma,\mu;\mu_1,\dots,\mu_{p'})\to
X_p^{m_1,\dots,m_p}(\Gamma,\mu;\lambda_1,\dots,\lambda_p).
$$
In view of (6), we obtain the formula
\begin{multline*}
\Phi\circ\Psi(e(nu,u_1,\dots,u_{p'})\\
{}=\sum_{\begin{subarray}{c}1\le\alpha_1\le m_1\\ \dots\\
1\le\alpha_p\le m_p\end{subarray}}
\Delta_{n,n-k}(u_1+\gamma_1u,\dots,u_{p'}+\gamma_{p'}u;
y_{\alpha_1,1}+\nu_1u,\dots,y_{\alpha_p,p}+\nu_pu)
e_{\alpha_1,\dots,\alpha_p}(nu).
\end{multline*}

    2. The function  $\theta(z)$ can be degenerated trigonometrically and rationally. In these cases, it
is replaced by the expressions $1-e^{2\pi
iz}$  and~$z$, respectively. Making these changes in the definitions
of $X_p^{m_1,\dots,m_p}(\Gamma,\mu;\lambda_1,\dots,\lambda_p)$ and $Y_{p'}(\Gamma,\mu;\mu_1,\dots,\mu_{p'})$, we obtain trigonometric and rational
degenerations of these algebras. Propositions 6 and 7 turn in this case into the constructions of
trigonometric and rational $R$-matrices if $\Theta_{n/k}(\Gamma)$ is replaced by the space of polynomials in the
variables $t_1,\dots,t_p$ of degree lower than~$n_j$ with respect to
the variable~$t_j$ ($1\le t\le p$), where
$t_j=e^{2\pi iz_j}$ in the trigonometric case and $t_j = z_j$ in the rational case.
\end{remarks}

\bigskip
\begin{flushleft}
L. D. Landau Institute of Theoretical Physics,\\
Russian Academy of Sciences
\end{flushleft}
\end{document}